\documentclass{article}
\usepackage[utf8]{inputenc}
\usepackage{lineno}
\usepackage{amssymb,mathrsfs, amsmath}
\usepackage{amsmath,amsthm}
\usepackage{amssymb}
\usepackage{latexsym}
\usepackage{amsfonts}
\usepackage{mathrsfs}
\usepackage{epsfig}
\usepackage{hyperref}
\usepackage{graphicx}
\usepackage{graphics}
\usepackage{pst-node}
\usepackage{stmaryrd}
\usepackage{cancel}
\usepackage {appendix}
\usepackage{titlefoot}
\usepackage{authblk}
\usepackage{tikz-cd}
\newtheorem{thm}{Theorem}[section]

\newtheorem{defn}{Definition}[section]
\newtheorem{prop}{Proposition}[section]
\newtheorem{lem}{Lemma}[section]

\newtheorem{rem}{Remark}[section]

\newtheorem{cor}{Corollary}[section]

\title{Cohomology and deformations of compatible Leibniz algebras}

\author[1]{R.B. Yadav}
\author[2]{Rinkila Bhutia}
\author[3]{Namita Behera}
\affil[1]{GITAM University, Bengaluru, Karnataka, 561203}
\affil[2,3]{Sikkim University, Gangtok, Sikkim, 737102, \textsc{India}}

\begin{document}

\maketitle
\begin{abstract} 
 In this paper we study a cohomology theory of compatible Leibniz algebra. We construct a bi-differential graded Lie algebra whose Maurer-Cartan elements characterize the  compatible Leibniz algebra structures. Using this, we study cohomology, infinitesimal deformations, Nijenhuis operator and their relation for compatible Leibniz algebras. Finally using cohomology of compatible Leibniz algebra with coefficients in an arbitrary representation we study the abelian extensions of compatible Leibniz algebra.
\end{abstract}

\vspace{0.3cm}

\vspace{0.3cm}

AMS: 17B56, 13D10, 17A30.


\unmarkedfntext{\newline\hspace*{2em} \textit{\thanks{} {\tt  rbyadav15@gmail.com, rbhutia@cus.ac.in,  nbehera@cus.ac.in}}}

\section{Introduction}
Leibniz algebra is a non-anti symmetric generalisation of Lie algebra. It was introduced and called D-algebra in papers by A. M. Bloch published in the 1960s to signify its relation with derivations. Later in 1993 J. L. Loday \cite{Loday} introduced the same structure and called it Leibniz algebra. Cohomology theory of Leibniz algebra with coefficients in a bimodule has been studied in \cite{LodayPiras}.

Algebraic deformation theory was introduced by Gerstenhaber for rings and algebra in a series of papers \cite{MG1}-\cite{MG5}. Subsequently algebraic deformation theory has been studied for different kind of algebras. To study deformation theory of any algebra, one needs a suitable cohomology, known as the deformation cohomology, which controls the deformation. 
In \cite{Bala}, D. Balavoine studies the formal deformation of algebras using the theory of Maurer-Cartan elements in a graded Lie algebra. In particular this approach is used to study the deformation of Leibniz algebra.

Here, we have defined a compatible Leibniz algebra to be a pair of Leibniz algebras such that the linear combination of their algebraic structure is also a Leibniz algebra. Recently cohomology and infinitesimal deformations of compatible Lie algebra and compatible associative algebra has been studied in \cite{CLieA} and \cite{CAA} respectively. Motivated by these works, in this paper we study the cohomology theory of compatible Leibniz algebra. Using the Balavoine bracket we define a graded Lie algebra whose Maurer-Cartan elements characterize the structure of compatible Leibniz algebras. We then study the cohomology of a compatible Leibniz algebra with coefficients in itself. This is then used to study infinitesimal deformation of compatible Leibniz algebra. We also establish the relation between Nijenhuis operator and the trivial infinitesimal deformation. Further we introduce the cohomology of compatible Leibniz algebra with coefficients in an arbitrary representation. Using this we study the abelian extensions of compatible Leibniz algebra.

This paper is organised as follows:
In section 2 we start with some basic concepts of Leibniz algebra. We then review Balavoine bracket, some results on cohomologies and the differential graded Lie algebra that controls the deformation of Leibniz algebra. In section 3 we define compatible Leibniz algebra and compatible bimodules. We then construct the graded Lie algebra whose Maurer-Cartan elements characterize compatible Leibniz algebra structure. In section 4 infinitesimal deformation of compatible Leibniz algebra is studied using cohomology of compatible Leibniz algebra with coefficients in itself. It is shown that equivalent infinitesimal deformations are in the same cohomology group. Then the notion of Nijenhuis operator on a compatible Leibniz algebra is studied and the correspondence between Nijenhuis operator and a trivial deformation is established. In section 5, cohomology of compatible Leibniz algebra with coefficients in an arbitrary representation is introduced. Finally in section 6, using the theory developed in section 5, abelian extension of compatible Leibniz algebra is studied. We end the paper by showing that the abelian extensions are classified by the second cohomology group.

Throughout the paper we consider the underlying field $K$ to be of characteristic $0$.

\section{Background}
\defn A Leibniz algebra is a vector space $L$ together with a $K$-linear operation $[.,.]: L\otimes L \to L$ such that
	$$[x,[y,z]]=[[x,y],z]+[y,[x,z]], ~\forall x,y, z \in L$$
	\defn A homomorphism between two Leibniz algebras $(L_1,[~]_1)$ and $(L_2,[~]_2)$ is a $K$-linear map $\phi:L_1\to L_2$ satisfying
$$\phi([x,y]_1)=[\phi(x),\phi(y)]_2. $$

\defn Let $(L,[~])$ be a Leibniz algebra. An $L$-bimodule is a vector space $M$ together with two $L$-actions
$m_l:L\otimes M\to M,~~ m_r: M\otimes L\to M$
such that for any $x,y\in L$ and $m\in M$ we have

$$m_l(x,m_l(y,m))=m_l([x,y],m)+m_l(y,m_l(x,m))$$
$$m_l(x,m_r(m,y))=m_r(m_l(x,m),y)+m_r(m,[x,y])$$	
$$m_r(m, [x,y])=m_r(m_r(m,x),y)+m_l(x, m_r(m,y)).$$

The following is a well established result.

\begin{prop} Let $(L,[~])$ be a Leibniz algebra and $(M,m_l,m_r)$ an $L$-bimodule. Then $L\oplus M$ is a Leibniz algebra with the Leibniz bracket defined as
$$[(x,u),(y, v)]_\ltimes=([x,y],m_l^1(x,v)+m_r^1(u,y))~~\forall~~x,y\in L~\text{and}~u,v\in M.$$
This is known as the semi-direct product.
\end{prop}

	\defn A permutation $\sigma \in S_n$ is called an $(i,n-i)$-shuffle if $\sigma(1)<\sigma(2)<...<\sigma(i)$ and $\sigma(i+1)<\sigma(i+2)<...<\sigma(n)$. If $i=0 ~\text{or}~ n$, we assume $\sigma=id$. $S_{(i, n-i)}$ denotes the set of all $(i,n-i)$-shuffles.
	
	\defn Let $(\mathfrak g=\oplus _{k\in \mathbb Z}\mathfrak g^k,[~], d)$ be a differential graded Lie algebra. A degree 1 element $x\in \mathfrak g^1$ is called a Maurer-Cartan element of $\mathfrak g$ if it satisfies
	$$dx+\frac{1}{2}[x,x]=0.$$
	
	\begin{thm} \cite{review} Let $(\mathfrak g=\oplus _{k\in \mathbb Z}\mathfrak g^k,[~])$ be a graded Lie algebra and $\mu \in \mathfrak{g^1}$ be a Maurer-Cartan element. Then the map
	$$d_\mu:\mathfrak g\to \mathfrak g,~d_{\mu}(u):=[\mu,u], ~\forall u \in \mathfrak g,$$
	 is a differential on $\mathfrak g$. \\
	 Further, for any $v \in \mathbb \mathfrak g^1$, the sum $\mu+v$ is a Maurer-Cartan element of the graded Lie algebra $(\mathfrak g=\oplus _{k\in \mathbb Z}\mathfrak g^k,[~])$ iff $v$ is a Maurer-Cartan element of the differential graded Lie algebra $(\mathfrak g=\oplus _{k\in \mathbb Z}\mathfrak g^k,[~], d_\mu).$
\end{thm}

	 \subsection {The Balavoine bracket}\cite{review}
	 Let $\mathfrak {g}$ be a vector space. We denote $ C^n(\mathfrak {g}, \mathfrak {g})=Hom(\otimes ^n\mathfrak {g}, \mathfrak {g})$ and set $\mathbb{C}^*(\mathfrak{g}, \mathfrak {g})=\oplus_{n\in \mathbb N}C^n(\mathfrak {g}, \mathfrak {g})$.   
  
  We assume the degree of an element in $C^n(\mathfrak {g}, \mathfrak {g})$ is $n-1$.
	 
	 For $P\in C^{p+1}(\mathfrak g, \mathfrak g), Q\in C^{q+1}(\mathfrak g, \mathfrak g)$ we define the \textbf{ Balavoine bracket} as
	 $$[P,Q]_B=P\circ Q-(-1)^{pq}Q\circ P$$
	where $P \circ Q \in  C^{p+q+1} $ is defined as \\
	$$(P\circ Q)(x_1,x_2,...,x_{p+q+1})=\sum_{k=1}^{p+1}(-1)^{(k-1)q}P\circ_k Q,$$
	 and 
  \begin{scriptsize}
  \begin{eqnarray*}
     && Po_kQ(x_1,x_2,...,x_{p+q+1})\\
	 &=&\sum_{\sigma \in S(k-1, q)}(-1)^{\sigma}P(x_{\sigma(1)},...,x_{\sigma(k-1)}, Q(x_{\sigma(k)},...,x_{\sigma(k+q-1)},x_{k+q}), x_{k+q+1},...,x_{p+q+1}).
  \end{eqnarray*}
  
	 \end{scriptsize}
	 \thm \label{I} The graded vector space $\mathbb C^*(\mathfrak g, \mathfrak g)$ equipped with the Balavoine bracket given above is a graded Lie algebra.
\\	  
In particular for $\pi \in  C^1(\mathfrak g, \mathfrak g)$, we have $[\pi,\pi]_B\in C^3(\mathfrak g, \mathfrak g) $ such that\\
$[\pi,\pi]_B=\pi \circ \pi -(-1)^{1.1}\pi \circ \pi= 2 \pi \circ \pi= 2\sum_{k=1}^{2}(-1)^{k-1}\pi \circ_k \pi= 2(\pi \circ_1 \pi - \pi \circ _2 \pi)$ \\
$\pi \circ_1 \pi(x,y,z)=\pi(\pi(x,y),z)$ and $\pi \circ_2 \pi(x,y,z)=\pi(x, \pi(y,z))-\pi(y, \pi(x,z))$	\\

Thus we have the following corollary. 

\cor \label{MC} $\pi$ defines a Leibniz algebra structure on $\mathfrak g$ iff $\pi$ is a Maurer-Cartan element of the graded Lie algebra $(\mathbb C^*(\mathfrak g, \mathfrak g),[~]_B).$

\thm \label{II} Let $(\mathfrak g, \pi)$ be a Leibniz algebra. Then $(\mathbb C^*(\mathfrak g, \mathfrak g), [~], d_\pi)$ becomes a differential graded Lie algebra (dgLa), where $d_\pi:=[\pi,.]_B$.\\
Further given $\pi' \in C^2(\mathfrak g, \mathfrak g)$, $\pi+\pi'$ defines a Leibniz algebra structure on $\mathfrak g$ iff $\pi'$ is a Maurer-Cartan element of the dgLa $(\mathbb C^*(\mathfrak g, \mathfrak g), [~], d_\pi)$.

\section{Compatible Leibniz algebra}

\defn A Compatible Leibniz algebra is a triple $(L,[~],\{~\})$ such that $(L,[~])$ and $(L,\{~\})$ are Leibniz algebras such that 	\begin{equation}
[x,\{y,z\}]+\{x,[y,z]\}=[\{x,y\},z]+\{[x,y],z\}+[y,\{x,z\}]+\{y,[x,z]\}    , ~\forall x,y, z \in L \label{CLA}
\end{equation}

\prop A triple $(L,[~],\{~\})$ is a compatible Leibniz algebra iff $(L,[~])$ and $(L,\{~\})$ are Leibniz algebras such that for any $k_1,k_2$ in $K$, the bilinear operation 
$$\llbracket x,y\rrbracket=k_1[x,y]+k_2\{x,y\}, ~\forall x,y \in L $$
defines a Leibniz algebra structure on L.
\begin{proof}
Let $(L,[~~],\{~~\})$ be a compatible Leibniz algebra. Then by definition itself $(L,[.,.])$ and $(L,\{.,.\})$ are Leibniz algebras. 
Further,
\begin{eqnarray}
&&\llbracket \llbracket x,y\rrbracket,z\rrbracket+\llbracket y,\llbracket x,z\rrbracket \rrbracket \notag\\
&=&\llbracket k_1[x,y]+k_2\{x,y\},z\rrbracket+\llbracket y, k_1[x,z]+k_2\{x,z\} \rrbracket \notag\\
&=&k_1[k_1[x,y]+k_2\{x,y\},z]+k_2\{k_1[x,y]+k_2\{x,y\},z\}+\notag\\
&&k_1[y, k_1[x,z]+k_2\{x,z\}]+k_2\{y, k_1[x,z]+k_2\{x,z\}\}\notag \\
&=&k_1k_1[[x,y],z]+k_1k_2[\{x,y\},z]+k_2k_1\{[x,y],z\}+k_2k_2\{\{x,y\},z\}+\notag\\
&&k_1k_1[y,[x,z]]+k_1k_2[y,\{x,z\}]+k_2k_1\{y,[x,z]\}+k_2k_2\{y,\{x,z\}\}\notag\\
&=&k_1^2([[x,y],z]+[y,[x,z]])+k_2^2(\{\{x,y\},z\}+\{y,\{x,z\}\})\notag\\
&&k_1k_2([\{x,y\},z]+\{[x,y],z\}+[y,\{x,z\}]+\{y,[x,z]\})\notag\\
&=&k_1^2[x,[y,z]]+k_2^2\{x,\{y,z\}\}+k_1k_2([x,\{y,z\}]+\{x,[y,z]\})\notag\\
&=&k_1(k_1[x,[y,z]]+k_2[x,\{y,z\}])+k_2(k_2\{x,\{y,z\}\}+k_1\{x,[y,z]\})\notag\\
	&=&k_1[x,k_1[y,z]+k_2\{y,z\}]+k_2\{x,k_2\{x,k_2{y,z}\}+k_1[y,z]\}\notag\\
	&=&k_1[x,\llbracket y,z \rrbracket]+k_2\{x,\llbracket y,z \rrbracket\}\notag\\
	&=& \llbracket x,\llbracket y,z\rrbracket \rrbracket.	
\end{eqnarray}
The converse is straight forward.
\end{proof}

\defn A homomorphism between two compatible Leibniz algebras $(L_1,[~]_1, \{~\}_1)$ and $(L_2,[~]_2, \{~\}_2)$ is a k-linear map $\phi:L_1\to L_2$ satisfying
$$\phi([x,y]_1)=[\phi(x),\phi(y)]_2~~\text{and}~~\phi(\{x,y\}_1)=\{\phi(x),\phi(y)\}_2. $$

\begin{defn} Let $(L,[~],\{~\})$ be a compatible Leibniz algebra. A compatible $L$-bimodule is a vector space $M$ together with four $L-$actions
 $$m_l^1:L\otimes M\to M,~~~~~~~~~~  m_r^1: M\otimes L\to M$$
$$m_l^2:L\otimes M\to M,~~~~~~~~~~  m_r^2: M\otimes L\to M$$
such that
\begin{itemize}
\item $(M,m_l^1, m_r^1)$ is a bimodule over $(L,[~~])$.
\item $(M,m_l^2, m_r^2)$ is a bimodule over $(L,\{~~\})$.
\item the following compatibilities hold for all $x,y \in L,~m\in M$
\scriptsize{
\begin{eqnarray*}
    LLM:~~ &&m_l^1(x,m_l^2(y,m))+m_l^2(x,m_l^1(y,m))\\
    &=&m_l^1(\{x,y\},m)+m_l^2([x,y],m)+m_l^1(y,m_l^2(x,m))+m_l^2(y,m_l^1(x,m))
\end{eqnarray*}
\begin{eqnarray*}
LML:~~&&m_l^1(x,m_r^2(m,y))+m_l^2(x,m_r^1(m,y))\\
&=&m_r^1(m_l^2(x,m),y)+m_r^2(m_l^1(x,m),y)+m_r^1(m, \{x,y\})+m_r^2(m, [x,y])
\end{eqnarray*}
\begin{eqnarray*}
MLL:~~&&m_r^1(m, \{x,y\})+m_r^2(m, [x,y])\\
&=&m_r^1(m_r^2(m,x),y)+m_r^2(m_r^1(m,x),y)+m_l^1(x, m_r^2(m,y))+m_l^2(x,m_r^1(m,y))
\end{eqnarray*}	
}	
\end{itemize}
\end{defn}
\it{Note}: Any compatible Leibniz algebra $(L,[~],\{~\})$ is a compatible $L-$bimodule in which $m_l^1=m_r^1=[~]$ and $m_l^2=m_r^2=\{~\}$.\\
Equivalently in terms of endomorphisms, we can define a compatible $L-$bimodule to be a vector space $M$ together with maps
$$\rho_1^L,\rho_1^R,\rho_2^L,\rho_2^R: L \to \text{End (M)}, ~\text{such that}$$
\begin{enumerate}
\item $(L, \rho_1^L,\rho_1^R)$ is an $L-$bimodule over $(L,[~~])$.
\item $(L, \rho_2^L,\rho_2^R)$ is an $L-$bimodule over $(L,\{~~\})$.
\item the following compatibilities hold for all $x,y \in L$

 $$\rho_1^L\{x,y\}+\rho_2^L[x,y]=[\rho_1^L(x),\rho_2^L(y)]+[\rho_2^L(x),\rho_1^L(y)]$$ 
 $$\rho_1^R\{x,y\}+\rho_2^R[x,y]=[\rho_1^L(x),\rho_2^R(y)]+[\rho_2^L(x),\rho_1^R(y)]$$
$$\rho_1^R(y)\rho_2^L(x)+\rho_2^R(y)\rho_1^L(x)=-\rho_1^R(y)\rho_2^R(x)-\rho_2^R(y)\rho_1^R(x).$$
\end{enumerate}

\subsection{Maurer-Cartan characterisation of Compatible Leibniz algebra}

\defn \cite{CLieA} Let ($\mathfrak g, [~~],\delta_1)$ and $(\mathfrak g, [~~],\delta_2)$ be two differential graded Lie algebras. We call $(\mathfrak g,[~~], \delta_1, \delta_2)$ a bi-differential graded Lie algebra (b-dgLa) if $\delta_1$ and $\delta_2$ satisfy
$$\delta_1\delta_2+\delta_2\delta_1=0.$$

\prop \cite{CLieA} Let $(\mathfrak g, [~~],\delta_1)$ and $(\mathfrak g, [~~],\delta_2)$ be two differential graded Lie algebras. Then $(\mathfrak g,[~~], \delta_1, \delta_2)$ is a bi-differential graded Lie algebra iff for any $k_1$ and $k_2$ $\in K$, $(\mathfrak g,[~~], \delta_{k_1k_2})$ is a differential graded Lie algebra where $\delta_{k_1k_2}=k_1\delta_1+k_2\delta_2$.

\defn Let $(\mathfrak g,[~~], \delta_1, \delta_2)$ be a b-dgLa. A pair $(\pi_1, \pi_2)\in \mathfrak g_1\oplus \mathfrak g_1$ is called a Maurer-Cartan element of the b-dgLa $(\mathfrak g,[~~], \delta_1, \delta_2)$ if $\pi_1$ and $\pi_2$ are Maurer-Cartan elements of the dgLas $(\mathfrak g,[~~], \delta_1)$ and $(\mathfrak g,[~~], \delta_2)$ respectively, and 
$$\delta _2\pi_1+\delta _1\pi_2+[\pi_1,\pi_2]=0. $$

\prop A pair $(\pi_1, \pi_2)\in \mathfrak g_1\oplus \mathfrak g_1$ is a Maurer-Cartan element of the b-dgLa $(\mathfrak g,[~~], \delta_1, \delta_2)$ iff for any $k_1,k_2\in K$, $k_1\pi_1+k_2\pi_2$ is a Maurer-Cartan element of the dgLa $(\mathfrak g,[~~], \delta_{k_1k_2})$.

\thm \label{MC bdgLa} Let $L$ be a vector space and $\pi_1, \pi_2 \in  C^2(L,L)$. Then $(L, \pi_1, \pi_2)$ is a compatible Leibniz algebra iff $(\pi_1, \pi_2)$ is a Maurer-Cartan element of the b-dgLa $(\mathbb C^*(L,L), [~~]_B, \delta_1=0, \delta_2=0)$.
\begin{proof}
$(L, \pi_1, \pi_2)$ is a compatible Leibniz algebra gives $(L, \pi_1)$ and $(L, \pi_2)$ are Leibniz algebras. Hence we get $[\pi_1, \pi_1]_B=[\pi_2, \pi_2]_B=0$.\\
Further $\forall x,y, z \in L$ we have the compatibility condition,
\begin{eqnarray}\label{CC}
&&\pi_1(x,\pi_2(y,z))+\pi_2(x,\pi_1(y,z))\notag\\
&=&\pi_1(\pi_2(x,y),z)+\pi_2(\pi_1(x,y),z)+\pi_1(y,\pi_2(x,z))+\pi_2(y,\pi_1(x,z))
\end{eqnarray}
We note that $[\pi_1, \pi_2]_B=\pi_1\circ \pi_2+\pi_2\circ \pi_1,$ where 
\begin{eqnarray*}
\pi_1\circ \pi_2(x,y,z)&=&(\pi_1\circ_1 \pi_2-\pi_1\circ_2 \pi_2)(x,y,z)\\
&=&\pi_1(\pi_2(x,y),z)-\pi_1(x,\pi_2(y,z))+\pi_1(y,\pi_2(x,z))  
\end{eqnarray*}
and 
\begin{eqnarray*}
\pi_2\circ \pi_1(x,y,z)&=&(\pi_2\circ_1 \pi_1-\pi_1\circ_2 \pi_1)(x,y,z) \\
&=&\pi_2(\pi_1(x,y),z)-\pi_2(x,\pi_1(y,z))+\pi_2(y,\pi_1(x,z)). 
\end{eqnarray*}
i.e.,
\begin{eqnarray*}
 [\pi_1, \pi_2]_B(x,y,z)  &=&\pi_1(\pi_2(x,y),z)-\pi_1(x,\pi_2(y,z))+\pi_1(y,\pi_2(x,z))\\
 &&+\pi_2(\pi_1(x,y),z)-\pi_2(x,\pi_1(y,z))+\pi_2(y,\pi_1(x,z)).    
\end{eqnarray*}
We thus have that
$[\pi_1, \pi_2]_B=0$ is equivalent to the compatibility condition \ref{CC}.

\end{proof}

\thm \cite{CLieA} \label {MC I}Let $(\pi_1,\pi_2)$ be a Maurer-Cartan element of the b-dgLa $(\mathfrak g,[~~], \delta_1, \delta_2)$.\\
Define $d_1:=\delta_1+[\pi_1,\_]$ and $d_2:=\delta_2+[\pi_2,\_]$. Then  $(\mathfrak g,[~~], d_1, d_2)$ is a b-dgLa.\\
Further for any $\Tilde{\pi}_1, \Tilde{\pi}_2 \in \mathfrak g_1$, $(\pi_1+\Tilde{\pi}_1, \pi_2+\Tilde{\pi}_2)$ is a Maurer-Cartan element of the b-dgLa $(\mathfrak g,[~~], \delta_1, \delta_2)$ iff $(\Tilde{\pi}_1,\Tilde{\pi}_2)$ is a Maurer-Cartan element of the b-dgLa $(\mathfrak g,[~~], d_1, d_2)$.\\

Let $(L,\pi_1, \pi_2)$ be a compatible Leibniz algebra. From theorems \ref{MC bdgLa} and \ref {MC I}, we conclude the following important results:

\thm $(\mathbb C^*(L,L),[~~], d_1,d_2)$ is a b-dgLa where  $d_1$ and $d_2$ are given by $d_1:=[\pi_1,\_]_B$ and $d_2:=[\pi_2,\_]_B$.
\thm  For any $\Tilde{\pi}_1, \Tilde{\pi}_2 \in C^2(L,L)$, $(L,\pi_1+\Tilde{\pi}_1, \pi_2+\Tilde{\pi}_2)$ is a compatible Leibniz algebra iff $(\pi_1+\Tilde{\pi}_1, \pi_2+\Tilde{\pi}_2)$ is a Maurer-Cartan element of the b-dgLa $(\mathbb C^*(L,L),[~~]_B, d_1, d_2)$.

\subsection{Cohomology of compatible Leibniz algebra}
Let $(L,[~~],\{~~\})$ be a compatible Leibniz algebra with $\pi_1(x,y)=[x,y]$ and $\pi_2(x,y)=\{x,y\}$.
By theorem \ref{MC bdgLa}, $(\pi_1,\pi_2)$ is a Maurer-Cartan element of the b-dgLa $(\mathbb C ^*(L,L), [~~]_B,0,0)$.

We define the cochains as follows:
 For $n\geq 1$, 
$$LC^n(L,L):=  \underbrace{\mathbb C^n(L,L) \oplus \mathbb C^n(L,L)...\oplus \mathbb C^n(L,L)}_{\text{n times} }$$
and $d^n:LC^n(L,L) \to LC^{n+1}(L,L)$ by
$$d^1f=([\pi_1,f]_B,[\pi_2,f]_B), ~\forall f\in LC^1(L,L)$$
$$d^n(f_1,f_2,..., f_n)=(-1)^{n-1}([\pi_1, f_1]_B,...,[\pi_2, f_{i-1}]_B+[\pi_1, f_i]_B,...,[\pi_2, f_n]_B),$$
$\text{where} (f_1,f_2,...,f_n)\in LC^n(L,L) ~\text{and}~ 2\leq i\leq n$.
 \\
  $d$ defined as above gives the following theorem. 
 \thm \label{d2} We have $d^{n+1}\circ d^n=0$.
  \begin{proof}
  We first note that since $(\pi_1,\pi_2)$ is a Maurer-Cartan element of the b-dgLa $(\mathbb C ^*(L,L), [~~]_B,0,0)$ we have
  $[\pi_1,\pi_1]_B=0,~[\pi_1,\pi_2]_B=0,~[\pi_2,\pi_2]_B=0$.\\
  For any $(f_1,f_2,\cdots,f_n)\in LC^n(L,L),~2\leq i \leq n$ we have
  
  \scriptsize{
  \begin{eqnarray}
   &&d^{n+1}d^n(f_1,f_2,\cdots,f_n)\notag\\
   &=&(-1)^{n-1}d^{n+1} ([\pi_1,f_1]_B,\cdots, [\pi_2,f_{i-1}]_B+[\pi_1,f_i]_B,\cdots,[\pi_2,f_n]_B)\notag\\
   &=&-([\pi_1,[\pi_1,f_1]_B]_B,[\pi_2,[\pi_1,f_1]_B]_B+[\pi_1,[\pi_2,f_1]_B]_B+[\pi_1,[\pi_1,f_2]_B]_B,\cdots\notag\\
   &&[\pi_2,[\pi_2,f_{i-2}]_B]_B+[\pi_2,[\pi_1,f_{i-1}]_B]_B+[\pi_1,[\pi_2,f_{i-1}]_B]_B+[\pi_1,[\pi_1,f_i]_B]_B,\cdots,\notag\\
   && [\pi_2,[\pi_2,f_{n-1}]_B]_B+[\pi_2,[\pi_1,f_{n}]_B]_B)+[\pi_1,[\pi_2,f_{n}]_B]_B, [\pi_2,[\pi_2,f_{n}]_B]_B)~~( 3\leq i\leq n-1)\notag\\
   &=&-(\frac{1}{2}[[\pi_1,\pi_1]_B,f_{1}]_B,[[\pi_1,\pi_2]_B,f_{1}]_B+\frac{1}{2}[[\pi_1,\pi_1]_B,f_{2}]_B,\cdots\notag\\
   &&\frac{1}{2}[[\pi_2,\pi_2]_B,f_{i-2}]_B+[[\pi_1,\pi_2]_B,f_{i-1}]_B+\frac{1}{2}[[\pi_1,\pi_1]_B,f_{i}]_B,\cdots,\notag\\
  && \frac{1}{2}[[\pi_2,\pi_2]_B,f_{n-1}]_B+[[\pi_1,\pi_2]_B,f_{n}]_B,\frac{1}{2}[[\pi_2,\pi_2]_B,f_{n}]_B)\notag\\
   &=&(0,0,\cdots,0).\notag
  \end{eqnarray}
 }
 	\end{proof}
Hence we have that $(LC^*(L,L))=(\oplus_{n\in \mathbb N}LC^n(L,L), d^*)$ is a cochain complex.

\defn Let $(L,[~~], \{~~\})$ be a compatible Leibniz algebra. The cohomology of the cochain complex $(LC^*(L,L), d^*)$ is called the cohomology of $(L,[~~], \{~~\})$. We denote the cohomology group by $H^n(L,L).$

\section{Infinitesimal deformations of compatible Leibniz algebras}
\defn Let $(L,[~~], \{~~\})$ be a compatible Leibniz algebra.
A \textit{formal one-parameter deformation} of  $L$ is a pair of $k[[t]]$-linear maps
$$\mu_t : L[[t]]\otimes L[[t]]\to L[[t]]~ \text{and}$$
$$m_t : L[[t]]\otimes L[[t]]\to L[[t]] ~\text{such that}:$$
\begin{itemize}
  \item[(a)] 
  $\mu_t(a,b)=\sum_{i=0}^{\infty}\mu_i(a,b) t^i$, ~~~
 $m_t(a,b)=\sum_{i=0}^{\infty}m_i(a,b) t^i$  
 
for all $a,b\in L$,  where $\mu_i, m_i:L\otimes L\to L$ are k-linear and $\mu_0(a,b,c)=[a,b]$ and $m_0(a,b)=\{a,b\}.$
 
      \item[(b)] For any $t$, $(L[[t]],\mu_t,m_t)$ is a compatible Leibniz algebra.
\end{itemize}
\defn Let $(L,[~~], \{~~\})$ be a compatible Leibniz algebra. Let $\mu_1, m_1\in C^2(L,L)$. Define
$$\mu_t(x,y)=[x,y]+t\mu_1(x,y),~~~~m_t(x,y)=\{x,y\}+t m_1(x,y),~~\forall x,y\in L.$$
If for any t, $(L,\mu_t, m_t)$ is a compatible Leibniz algebra, we say that $(L,\mu_t,m_t)$ defines an infinitesimal deformation of $(L,[~~], \{~~\})$.\\
We also say that $(\mu_1,m_1)$ generates an infinitesimal deformation of $(L,[~~], \{~~\})$.\\
For convenience we write $[x,y]=\mu_0(x,y)$ and $\{x,y\}=m_0(x,y)$.\\
By \ref{MC bdgLa} we have that $(L,\mu_t, m_t)$ is a compatible Leibniz algebra if and only if 
$(\mu_t, m_t)$ is a Maurer-Cartan element of $(C^*(L,L),[~~]_B,0,0)$.  ``$(\mu_t, m_t)$ is a Maurer-Cartan element " is equivalent to following condition
\begin{equation}\label{rnr1}
  [\mu_t,\mu_t]_B=0~~[m_t,m_t]_B=0~~[\mu_t,m_t]_B=0  
\end{equation}
Condition \ref{rnr1} is equivalent to following conditions
 \begin{equation*}
     [\mu_0,\mu_0]_B=0,~[\mu_0,\mu_1]_B=0, ~[\mu_1,\mu_1]_B=0
 \end{equation*}
\begin{equation*}
    [m_0,m_0]_B=0,~[m_0,m_1]_B=0,~[m_1,m_1]_B=0
\end{equation*}
\begin{equation*}
  [\mu_0,m_0]_B=0,~[\mu_0,m_1]_B+[\mu_1,m_0]_B=0,~[\mu_1,m_1]_B=0  
\end{equation*}
	
Reordering the terms and excluding the trivial equations we get that
$(L,\mu_t,m_t)$ defines an infinitesimal deformation of $(L,[~~], \{~~\})$ iff
$$[\mu_0,\mu_1]_B=0,~~~~[m_0,m_1]_B=0~,~~~~[\mu_0,m_1]_B+[\mu_1,m_0]_B~=0$$
$$[\mu_1,\mu_1]_B=0,~~~~[m_1,m_1]_B=0,~~~~[\mu_1,m_1]_B=0~.$$

Note that the first line above implies $d^2(\mu_1,m_1)=0$ i.e $(\mu_1,m_1)$ is a 2-cocycle and the second line implies that $(L,\mu_1,m_1)$ is a compatible Leibniz algebra.

Hence we have the following theorem.
\thm  Let $(L,[~~], \{~~\})$ be a compatible Leibniz algebra. If $(\mu_1,m_1)\in LC^2(L,L)$ generates an infinitesimal deformation then $(\mu_1,m_1)$ is a cocycle.

\defn Two infinitesimal deformations $(L,\mu_t,m_t)$ and $(L,\mu_t',m_t')$ are said to be equivalent if there exists a linear bijection $N:L\to L$ such that
$$Id+tN:(L,\mu_t,m_t) \to (L,\mu_t',m_t')$$
is a compatible Leibniz algebra homomorphism.

$Id+tN$ being a compatible Leibniz algebra homomorphism implies
\begin{enumerate}
	\item $[x,y]=[x,y]'$
	\item $\mu_1(x,y)-\mu_1'(x,y)=[x,N(y)]+[N(x),y]-N[x,y]$
	\item $N\mu_1(x,y)=\mu_1'(x,N(y))+\mu_1'(N(x),y)+[N(x),N(y)]$
	\item $\mu_1'(N(x),N(y))=0$
	\item $\{x,y\}=\{x,y\}'$
	\item $m_1(x,y)-m_1'(x,y)=\{x,N(y)\}+\{N(x),y\}-N\{x,y\}$
	\item $N m_1(x,y)=m_1'(x,N(y))+m_1'(N(x),y)+\{N(x),N(y)\}$
	\item $m_1'(N(x),N(y))=0$
	\end{enumerate}

$2$ and $6$ gives
\begin{eqnarray}
	(\mu_1-\mu_1', m_1-m_1')(x,y)&=&([x,N(y)]+[N(x),y]-N[x,y], \{x,N(y)\}\notag\\
 &&+\{N(x),y\}-N\{x,y\})\notag\\
	&=&([\mu_0,N]_B,[m_0,N]_B)(x,y)\notag\\
	&=&d^1N(x,y). \notag
	\end{eqnarray}
Thus we have the following theorem\\
\thm If two infinitesimal deformations $(L,\mu_t,m_t)$ and $(L,\mu_t',m_t')$ of a compatible Leibniz algebra $(L,\mu_0,m_0)$ are equivalent then, $(\mu_1, m_1)$ and $(\mu_1',m_1')$ are in the same cohomology class.

\defn Let $(L,[~~])$ be a Leibniz algebra. A linear map $N:L\to L$ is said to be a Nijenhuis operator on $L$ if 
$$N([x,N(y)]+[N(x),y]-N[x,y])=[N(x),N(y)]~~\forall x,y \in L.$$
We define linear $[~~]_N: L\otimes L \to L$ as 
$$[x,y]_N=[x,N(y)]+[N(x),y]-N[x,y]$$
$T_{[~~]_N}: L\otimes L\to L$ denotes the Nijenhuis torsion of $N$ defined as\\
$$T_{[~~]_N}(x,y)=N([x,y]_N)-[N(x), N(y)], ~~\forall x,y \in L.$$
When $N$ is a Nijenhuis operator we get that $T_{[~~]_N}=0$.

\prop \label{Nijen} If $N:L\to L$ is a Nijenhuis operator on Leibniz algebra $(L,[~~])$, then $(L,[~~]_N)$ is also a Leibniz algebra. Further $N$ is a Leibniz algebra homomorphism from $(L,[~~]_N)$ to $(L,[~~])$.
Furthermore $(L,[~~],[~~]_N)$ forms a compatible Leibniz algebra.

\begin{proof} 
 $$ [x,[y,z]_N]_N=[[x,y]_N,z]_N+[y,[x,z]_N]_N, ~\forall x,y, z \in L$$
Further, $N([x,y]_N)=[N(x),N(y)]$ follows from the definition of Nijenhuis operator and $[~~]_N$.

To show $(L,[~~],[~~]_N)$ is a compatible Leibniz algebra we first note that $\pi_N=[\pi, N]_B$. For any $k_1$ and $k_2\in K$,
\begin{eqnarray*}
[k_1\pi+k_2\pi_N,k_1\pi+k_2\pi_N]_B&=&k_1k_2([\pi,\pi_N]_B+[\pi_N,\pi]_B)\\
&=&2k_1k_2[\pi,\pi_N]_B\\
&=&2k_1k_2[\pi,[\pi,N]_B]_B\\
&=&0.    
\end{eqnarray*}
\end{proof}

\defn  Let $(L,[~~],\{~~\})$ be a compatible Leibniz algebra. A linear map $N:L\to L$ is said to be a Nijenhuis operator on $(L,[~~],\{~~\})$ if $N$ is a Nijenhuis operator on the Leibniz algebras $(L,[~~])$ and $(L,\{~~\})$.

\prop Let $(L,[~~],\{~~\})$ be a compatible Leibniz algebra. A linear map $N:L\to L$ is a Nijenhuis operator on $(L,[~~],\{~~\})$ iff for any $k_1,k_2$ in $K$, $N$ is a Nijenhuis operator on the Leibniz algebra $(L,\llbracket~~\rrbracket)$, where
$\llbracket x,y \rrbracket=k_1[x,y]+k_2\{x,y\}, ~\forall x,y \in L$.
\proof We have  \begin{eqnarray}
T_{\llbracket~~\rrbracket_N}(x,y)&=&N(\llbracket x,y \rrbracket_N)-\llbracket N(x),N(y) \rrbracket\notag\\
	&=&N(k_1[x,Ny]+k_2\{x,Ny\}+k_1[Nx,y]+k_2\{Nx,y\})\notag\\
 &&-k_1[N(x), N(y)]-k_2\{N(x), N(y)\}\notag\\
	&=&k_1(N([x,y]_N)-[N(x), N(y)])+k_2(N(\{x,y\}_N)-\{N(x), N(y)\})\notag\\
	&=&k_1T_{[~~]_N}(x,y)+k_2T_{\{~~\}_N}(x,y)\notag\\
\end{eqnarray}
Hence we have, 
$$T_{\llbracket~~\rrbracket_N}=0~~\text{iff}~~~T_{[~~]_N}=T_{\{~~\}_N}=0.$$

\prop \label{Nijen1}  Let $(L,[~~],\{~~\})$ be a compatible Leibniz algebra and $N:L\to L$ be a Nijenhuis operator on $(L,[~~],\{~~\})$. Then $(L,[~~]_N,\{~~\}_N)$ is also a compatible Leibniz algebra and $N$ is a compatible Leibniz algebra homomorphism from $(L,[~~]_N,\{~~\}_N)$ to $(L,[~~],\{~~\})$.
\proof Let $N:L\to L$ be a Nijenhuis operator on $(L,[~~],\{~~\})$. then by the previous theorem  $N$ is a Nijenhuis operator on the Leibniz algebra $(L,\llbracket~~\rrbracket)$ for any $k_1,k_2$ in $K$.\\
Using result \ref{Nijen} we get that $(L,\llbracket~~\rrbracket_N)$ is a Leibniz algebra and $N$ is a Leibniz algebra homomorphism from $(L,\llbracket~~\rrbracket_N)$ to $(L,\llbracket~~\rrbracket)$.\\
Hence we have that $(L,[~~]_N,\{~~\}_N)$ is a compatible Leibniz algebra. And we also get that
$N$ is a compatible Leibniz algebra homomorphism from $(L,[~~]_N,\{~~\}_N)$ to $(L,[~~],\{~~\})$.

\defn An infinitesimal deformation $(L,\mu_t,m_t)$ of compatible Leibniz algebra $(L,\mu_0,m_0)$ generated by $(\mu_1,m_1)$ is trivial if there exists linear $N: L\to L$ such that $Id+tN:(L,\mu_t,m_t) \to (L,\mu_0,m_0)$ is a compatible Leibniz algebra homomorphism.

$Id+tN$ is a compatible Leibniz algebra homomorphism iff
\begin{enumerate}
\item $\mu_1(x,y)=[x,N(y)]+[N(x),y]-N[x,y]$
\item $m_1(x,y)=\{x,N(y)\}+\{N(x),y\}-N\{x,y\}$
\item $N\mu_1(x,y)=[N(x),N(y)]$
\item $N m_1(x,y)=\{N(x),N(y)\}$
\end{enumerate}
$1$ and $3$ gives that $N$ is a Nijenhuis operator on $(L,\mu_0)$. $2$ and $4$ gives that $N$ is a Nijenhuis operator on $(L,m_0)$.

Thus we have the following theorem.

\thm  \label{NO} A trivial infinitesimal deformation of a compatible Leibniz algebra gives rise to a Nijenhuis operator.


\thm A Nijenhuis operator on a compatible Leibniz algebra $(L,[~~],\{~~\})$ gives rise to a trivial deformation.
\proof Let $N$ be a  Nijenhuis operator on a compatible Leibniz algebra $(L,[~~],\{~~\})$. Take
$$\mu_1(x,y)=[x,N(y)]+[N(x),y]-N[x,y]$$
$$m_1(x,y)=\{x,N(y)\}+\{N(x),y\}-N\{x,y\}$$
for any $x,y\in L$. Then
\begin{eqnarray}
d^1N(x,y)&=&([\mu_0,N]_B,[m_0,N]_B)(x,y)\notag\\
&=&([x,N(y)]+[N(x),y]-N[x,y], \{x,N(y)\}+\{N(x),y\}-N\{x,y\})\notag\\
	&=&(\mu_1(x,y),m_1(x,y)).\notag
\end{eqnarray}
i.e., $(\mu_1,m_1)$ is a 2-cocycle.\\
Further since $N$ is a Nijenhuis operator on $(L,[~~],\{~~\})$, and $\mu_1=[~~]_N$ and $m_1=\{~~\}_N$, by proposition (\ref{Nijen1}) we get that $(L, [~~]_N,\{~~\}_N)$ is a compatible Leibniz algebra.\\
These two statements implies that $(\mu_1,m_1)$ give rise to an infinitesimal deformation of $L$.
Showing the deformation is trivial is straightforward.

\section{Cohomologies of compatible Leibniz algebras with coefficients in arbitrary representation}

For vector spaces $g_1$ and $g_2$, we define $g^{l,k}$ to be the direct sum of tensor products of $g_1$ and $g_2$, where $g_1$ is repeated $l$ times and $g_2$ is repeated $k$ times. For example $g^{1,1}=(g_1\otimes g_2) \oplus(g_2\otimes g_1) $ and $g^{2,1}=(g_1\otimes g_1\otimes g_2)\oplus(g_1\otimes g_2\otimes g_1)\oplus(g_2\otimes g_1\otimes g_1) $. Thus $\otimes ^n(g_1\oplus g_2)\equiv \oplus_{l+k=n}g^{l,k}$.\\ 
For any linear map $f:g_{i_1}\otimes g_{i_2}\cdots \otimes g_{i_n}\to g_j,~~\text{where} i_1,i_2,...,i_n,j\in\{1,2\}$ we define $\hat{f}\in C^n(g_1\oplus g_2,g_1\oplus g_2 )$ as
\[
    \hat{f}= 
\begin{cases}
    f,& \text{on } g_{i_1}\otimes g_{i_2}\cdots \otimes g_{i_n}\\
    0,              & \text{otherwise}
\end{cases}
\]
$\hat{f}$ is called a lift of $f$.\\
In particular, for the linear maps we encountered in the previous sections:
$$\pi:L\otimes L\to L,~~m_l:L\otimes M \to M,~~m_r:M\otimes L \to M$$
we get lifts
$$\hat{\pi}:(L\oplus M)^2\to L\oplus M,~\text{defined as}~ \hat{\pi}((x_1,v_1),(x_2,v_2))=(\pi(x_1,x_2),0) $$
$$\hat{m}_l:(L\oplus M)^2\to L\oplus M,~\text{defined as} ~\hat{m}_l((x_1,v_1),(x_2,v_2))=(0,m_l(x_1,v_2))$$
$$\hat{m}_r:(L\oplus M)^2\to L\oplus M,~\text{defined as} ~\hat{m}_r((x_1,v_1),(x_2,v_2))=(0,m_r(v_1,x_2))$$

By property  of the Hom-functor we get
$$C^n(g_1\oplus g_2, g_1\oplus g_2)\equiv \sum_{l+k=n}C^n(g^{l,k},g_1)\oplus \sum_{l+k=n}C^n(g^{l,k},g_2). $$

\defn A linear map $f\in Hom(\otimes ^n(g_1\oplus g_2), (g_1\oplus g_2) )$ has bidegree $l|k$ if
\begin{enumerate}
\item $l+k+1=n$
\item if $X\in g^{l+1,k}$ then $f(X)\in g_1$
\item if $X\in g^{l,k+1}$ then $f(X)\in g_2$
\item $f(X)=0$ in all other cases.
\end{enumerate}
We use notation $\Vert f \Vert=l|k$. We say that 
$f$ is homogeneous if $f$ has a bidegree.

Considering examples above, we have $\Vert \hat{\pi}\Vert=\Vert \hat{m}_l\Vert=\Vert \hat{m}_r\Vert=1|0$.\\

In the next three lemmas, we consider a few standard results regarding bidegrees \cite{rong}. 

\lem \label{L1} If $f_1,f_2,\cdots f_k \in C^n(g_1\oplus g_2, g_1\oplus g_2)$ be homogeneous linear maps and the bidegrees of $f_i$ are different. Then $f_1+f_2+...+f_k=0$ iff $f_1=f_2=\cdots=f_k=0$.

\lem \label{L2} If $\Vert f\Vert =-1/l ~(l/-1)$ and  $\Vert g\Vert =-1/k ~(k/-1)$ then $[f,g]_B=0$.

\lem \label{L3} $f\in C^n(g_1\oplus g_2, g_1\oplus g_2)$ and $g\in C^m(g_1\oplus g_2, g_1\oplus g_2)$ be homogeneous linear maps with bidegrees $l_f|k_f$ and $l_g|k_g$ respectively. Then $[f,g]_B$ is a linear map of bidegree $l_f+l_g|k_f+k_g$.
\thm \cite{review} Let $(L,\pi=[~~])$ be a Leibnz algebra. $(V,m_l, m_r)$ is a representation of $L$ iff  $\hat{m}_l+\hat{m}_r$ is a Maurer-Cartan element of the dgLA $(C^*(L\oplus V, L\oplus V),[~~]_B, \partial_{\hat{\pi}}=[\hat{\pi}, .]_B)$.\\

\cor \label{MC1} If $(V,m_l, m_r)$ is a representation of $(L,\pi)$, then $[\hat{\pi}+\hat{m}_l+\hat{m}_r,\hat{\pi}+\hat{m}_l+\hat{m}_r]_B=0$. 

Let $(V,m_l,m_r)$ be a representation of the Leibniz algebra $(L,\pi)$.\\
We denote the set of all homogeneous elements of degree $p|q$ of $C^{p+q+1}(L\oplus V, L\oplus V)$ by $C^{p|q}(L\oplus V, L\oplus V )$.
We define the set of $n$-cochains as
$$C^n(L,V):=C^{n|-1}(L\oplus V, L\oplus V)\cong Hom(\otimes ^n L,V) \text{ using the lift map}$$

and coboundary operator $d^n_{\pi+m_l+m_r}:C^n(L,V)\to C^{n+1}(L,V)$ as
$$d^n_{\pi+m_l+m_r}f:=(-1)^{n-1}[\hat{\pi}+\hat{m}_l+\hat{m}_r,\hat{f}]_B,~\forall f\in C^n(L,V).$$
Note that since $\hat{\pi}+\hat{m}_l+\hat{m}_r\in C^{1|0}$ and $\hat{f}\in C^{n|-1}$, Lemma \ref{L3} gives us that $[\hat{\pi}+\hat{m}_l+\hat{m}_r,\hat{f}]_B\in C^{n+1|-1}$.\\

Further note that
$d^{n+1}d^nf=-[\hat{\pi}+\hat{m}_l+\hat{m}_r,[\hat{\pi}+\hat{m}_l+\hat{m}_r,\hat{f}]_B]_B=0$ by the graded Jacobi identity.\\
Thus we have a well defined cochain complex $(C^*(L,V), d^*_{\hat{\pi}+\hat{m}_l+\hat{m}_r})$.
~~~~~~~~~~~~~~~~~~~~~~~~~~~~~~~~~~~~~~~~~~~~~~~~~~~~~~~~~~~~~~~~~~~~~~~~~~~~~~~~~~~~~~~~~~~~
\thm Let $(L,\pi_1=[~~],\pi_2=\{~~\})$ be a compatible Leibniz algebra and $(M,m^1_l,m^1_r,m^2_l,m^2_r )$ a representation of $L$. Then $(\hat{\pi}_1+\hat{m}^1_l+\hat{m}^1_r, \hat{\pi}_2+\hat{m}^2_l+\hat{m}^2_r)$ is a Maurer-Cartan element of the bi-differential graded Lie Algebra $(C^*(L\oplus M, L\oplus M), [~~]_B, 0, 0)$ i.e 
\begin{equation}\label{I}
[\hat{\pi}_1+\hat{m}^1_l+\hat{m}^1_r, \hat{\pi}_1+\hat{m}^1_l+\hat{m}^1_r]_B=0,
\end{equation}
\begin{equation}\label{2}
[\hat{\pi}_1+\hat{m}^1_l+\hat{m}^1_r, \hat{\pi}_2+\hat{m}^2_l+\hat{m}^2_r]_B=0,
\end{equation}
\begin{equation}\label{3}
[\hat{\pi}_2+\hat{m}^2_l+\hat{m}^2_r, \hat{\pi}_2+\hat{m}^2_l+\hat{m}^2_r]_B=0
\end{equation}
\proof Since  $(M,m^1_l,m^1_r)$ is a representation of the Leibniz algebra $(L,\pi_1=[~~])$, by corollary \ref{MC1}
equation \ref{I} holds. 
Likewise $(M,m^2_l,m^2_r)$ is a representation of the Leibniz algebra $(L,\pi_2=\{~~\})$, by corollary~\ref{MC1}
equation \ref{3} holds.\\
For $x_1,x_2,x_3 \in L$, $v_1, v_2, v_3 \in V$
\begin{eqnarray}
&&[\hat{\pi_1}+\hat{m}^1_l+\hat{m}^1_r, \hat{\pi_2}+\hat{m}^2_l+\hat{m}^2_r]_B(x_1,v_1),(x_2,v_2),(x_3,v_3)\notag\\
	&=&(\pi_1(\pi_2(x_1,x_2),x_3), m^1_l(\pi_2(x_1,x_2),u_3)+m^1_r(m^2_l(x_1,u_2)+(m^2_r(u_1,x_2),x_3))\notag\\
 &+&(-\pi_1(x_1,\pi_2(x_2,x_3),-m^1_l(x_1,m^2_l(x_2,u_3)+m^2_r(u_2,x_3))-m^1_r(u_1,\pi_2(x_2,x_3)))\notag\\
 &+&(\pi_1(x_2,\pi_2(x_1,x_3)),m^1_l(x_2,m^2_l(x_1,u_3)+m^2_r(u_1,x_3))+m^1_r(u_2,\pi_2(x_1,x_3)))\notag\\
  &+&(\pi_2(\pi_1(x_1,x_2),x_3),m^2_l(\pi_1(x_1,x_2),u_3)+m^2_r(m^1_l(x_1,u_2)+m^1_r(u_1,x_2),x_3)))\notag\\
   &+&(-\pi_2(x_1,\pi_1(x_2,x_3)),-m^2_l(x_1,m^1_l(x_2,u_3)+m^1_r(u_2,x_3))-m^2_r(u_1,\pi_1(x_2,x_3)))\notag\\
    &+&((\pi_2(x_2,\pi_1(x_1,x_3)),m^2_l(x_2,m^1_l(x_1,u_3)+m^1_r(u_1,x_3))+m^2_r(u_2,\pi_1(x_1,x_3)))\notag\\
    &=&0.
\end{eqnarray}
We get the above by the compatibility conditions \ref{CLA}, $LLM$, $LML$ and $MLL$.\\

Note that the  coboundary operator for $(L,\pi_1)$ with coefficients in $(M,m^1_l,m^1_r)$ and for $(L,\pi_2)$ with coefficients in $(M,m^2_l,m^2_r)$ are respectively given by

$$d^n_{\pi^1+m^1_l+m^1_r}f:=(-1)^{n-1}[\hat{\pi}_1+\hat{m}^1_l+\hat{m}^1_r,\hat{f}]_B,~\text{and}$$

$$d^n_{\pi^2+m^2_l+m^2_r}f:=(-1)^{n-1}[\hat{\pi}_2+\hat{m}^2_l+\hat{m}^2_r,\hat{f}]_B,~\forall f\in C^n(L,M)$$
By the graded Jacobi identity it can be shown that the three  conditions \ref{I} , \ref{2}, \ref{3} implies
\begin{eqnarray}
    d^{n+1}_{\hat{\pi}_1+\hat{m}^1_l+\hat{m}^1_r}d^n_{\hat{\pi}_1+\hat{m}^1_l+\hat{m}^1_r}=0\notag\\
    d^{n+1}_{\hat{\pi}_2+\hat{m}^2_l+\hat{m}^2_r}d^n_{\hat{\pi}_2+\hat{m}^2_l+\hat{m}^2_r}=0\notag\\
    d^{n+1}_{\hat{\pi}_1+\hat{m}^1_l+\hat{m}^1_r}d^n_{\hat{\pi}_2+\hat{m}^2_l+\hat{m}^2_r}+d^{n+1}_{\hat{\pi}_2+\hat{m}^2_l+\hat{m}^2_r}d^n_{\hat{\pi}_1+\hat{m}^1_l+\hat{m}^1_r}=0\label{*}
\end{eqnarray}
For $n\geq1 $ we define the space of n-cochains $ LC^n(L,M)$ as
$$LC^n(L,M)=C^n(L,M)\oplus C^n(L,M)\oplus \cdots \oplus C^n(L,M)~~ \}\text{n-copies} $$

and coboundary for $n\geq 1$, $\partial^n:LC^n \to LC^{n+1}$ as
$$\partial ^1 f= (d_{\hat{\pi}_1+\hat{m}^1_l+\hat{m}^1_r}f, d_{\hat{\pi}_2+\hat{m}^2_l+\hat{m}^2_r}f)~~\forall f\in Hom(L,M)$$
and for $2\leq i\leq n$ and $(f_1,f_2,\cdots, f_n)\in LC^n(L,V)$,
\begin{eqnarray*}
    &&\partial ^n(f_1,f_2,\cdots, f_n)\notag\\
    &= &(d^n_{\hat{\pi}_1+\hat{m}^1_l+\hat{m}^1_r}f_1,\cdots, d^n_{\hat{\pi}_2+\hat{m}^2_l+\hat{m}^2_r}f_{i-1}+d^n_{\hat{\pi}_1+\hat{m}^1_l+\hat{m}^1_r}f_i,\cdots, d_{\hat{\pi}_2+\hat{m}^2_l+\hat{m}^2_r}f_n).
\end{eqnarray*}
Using \ref{*} or just like in Theorem \ref{d2} it can be shown that $\partial ^2=0$.

\defn Let  $(M,m^1_l,m^1_r, m^2_l,m^2_r)$ be a representation of  a  compatible Leibniz algebra $(L,\pi_1,\pi_2)$. The cohomology of the cochain complex $(\oplus_{n=1}^\infty LC^n(L,L),\partial)$ is called the cohomology of $(L,\pi_1,\pi_2)$ with coefficient in the representation $(M,m^1_l,m^1_r, m^2_l,m^2_r)$.
The corresponding $n^{th}$cohomology group is denoted by $\mathbb H^n(L,M)$.

\section{Abelian extension of compatible Leibniz algebras}
\defn Let $(L,\pi_1,\pi_2)$ and $(V,\mu_1,\mu_2)$ be two compatible Leibniz algebras. An extension of $(L,\pi_1,\pi_2)$ by $(V,\mu_1,\mu_2)$ is a short exact sequence of compatible Leibniz algebra morphisms
\begin{equation}
0\to V\overset{i}\to \tilde {L}\overset{\rho}\to L\to 0,
\end{equation}
where $(\tilde L,[~~]_{\tilde L}, \{~\}_{\tilde L})$ is a compatible Leibniz algebra.

 A section of the extension $(\tilde L,\tilde \pi_1, \tilde \pi_2)$ of $(L,\pi_1,\pi_2)$ by $(V,\mu_1,\mu_2)$ is a linear map $s :L\to \tilde L$ such that $\rho \circ s =Id_L$ where $Id_L$ is the identity morphism on $L$.
\\
$\tilde L$ is called an abelian extension of $L$ by $V$, if the compatible Leibniz algebra structure on $V$ is trivial i.e., $\mu_1(x,y)=\mu_2(x,y)=0, ~\forall x,y \in V$.


\begin{defn}
Two  abelian extensions, $(L_1,[~~]_1, \{~\}_1)$ and $(L_2,[~~]_2, \{~\}_2)$ of $(L,\pi_1,\pi_2)$ by $(V,0,0)$ given by exact sequences 
$0\to V\overset{i}\to L_1\overset{\rho}\to L\to 0$ and $0\to V\overset{i'}\to L_2\overset{\rho'}\to L\to 0$
are equivalent if there exists an isomorphism $F: L_1\to  L_2$ such that the following diagram commutes

\begin{equation}\label{commdia}
\begin{tikzcd}
0\arrow{r} &  V \arrow{r}{i} \arrow{d}{Id_v}  &  L_1 \arrow{r}{\rho} \arrow{d}{F}  &  L \arrow{r}\arrow{d}{Id_L} & 0\\
0\arrow{r} &  V \arrow{r}{i'}                          &  L_2 \arrow{r}{\rho'}                             &L\arrow{r}                         & 0.
\end{tikzcd}
\end{equation}
\end{defn}

Consider an abelian extension $\tilde L$ of $L$ by $V$ and section $s:L\to \tilde L$ For $x,y\in L,~u\in V$ define linear maps \\
$\theta: L\otimes L \to V$, $\tilde{\theta}: L\otimes L \to V$, $m_l^1,~m_l^2:L\otimes V\to V$ and $m_r^1,  m_r^2 : V\otimes L\to V$
as
$$\theta(x,y)=[s(x),s(y)]_{\tilde L}-s[x,y]$$
$$\tilde{\theta}(x,y)=\{s(x),s(y)\}_{\tilde L}-s\{x,y\}$$
$$m_l^1(x,u)=[s(x),u]_{\tilde L},~~m_r^1(u,x)=[u,s(x)]_{\tilde L}$$
$$m_l^2(x,u)=\{s(x),u\}_{\tilde L}~~m_r^2(u,x)=\{u,s(x)\}_{\tilde L}.$$

\thm \label{Aex1} With the notations introduced above, $(V,m_l^1,m_r^1, m_l^2,m_r^2)$ is a representation of the compatible Leibniz algebra $(L, \pi_1,\pi_2)$. Further this representation is independent of the choice of sections $s$.
\proof  For $x,y\in L$ and $u\in V$, since $\theta(x,y) \in V$, we have
\begin{eqnarray}
0&=&[\theta(x,y),u]_{\tilde L}\notag\\
&=&[[s(x),s(y)]_{\tilde L},u]_{\tilde L}-[s[x,y],u]_{\tilde L}\notag\\
&=&[s(x),[s(y),u]_{\tilde L}]_{\tilde L}-[s(y),[s(x),u]_{\tilde L}]_{\tilde L}-[s[x,y],u]_{\tilde L}\notag\\
\end{eqnarray}
i.e., we have $m_l^1([x,y],u)=m_l^1(x,m_l^1(y,u))-m_l^1(y,m_l^1(x,u))$.\\
Similarly, we can show 
$$m_r^1(u,[x,y])=m_l^1(x,m_r^1(u,y))-m_r^1(m_l^1(x,u),y) \text{and}$$
$$m_r^1(u,[x,y])=m_r^1(m_r^1(u,x),y)+m_l^1(x,m_r^1(u,y)).$$
Thus $(V, m_l^1,m_r^1)$ is a representation of $(L,[~]).$\\
Similarly we can prove that $(V, m_l^2,m_r^2)$ is a representation of $(L,\{~\})$ and the compatibility conditions $LLM,~LML$ and $MLL$ hold.\\

Suppose there is another section $s':L\to \tilde L$ w.r.t which $(V,m_{l}^{1'},m_{r}^{1'}, m_{l}^{2'},m_{r}^{2'})$ is the representation of $(L,\pi_1,\pi_2)$.\\
Then for any $x\in L$ and $u\in V$
\begin{eqnarray}
    m_l^1(x,u)-m_l^{1'}(x,u)&=&[s(x),u]_{\tilde L}-[s'(x),u]_{\tilde L}\notag\\
    &=&[s(x)-s'(x),u]_{\tilde L}\notag\\
    &=&0,~~~ (\text{since } s(x)-s'(x)\in V) \notag
 \end{eqnarray}
Hence $m_l^1=m_l^{1'}$. \\
Similarly, it can be shown that $m_r^1= m_r^{1'},~m_l^2=m_l^{2'}, ~m_r^2=m_r^{2'}$.

\thm \label{Aex2} Two equivalent abelian extensions, $(L_1,[~~]_1, \{~\}_1)$ and $(L_2,[~~]_2, \{~\}_2)$ of $(L,\pi_1,\pi_2)$ by $(V,0,0)$ give rise to the same representation of $(L,\pi_1,\pi_2)$ as in the previous theorem.
\proof Let $s_1$ and $s_2$ be sections of $(L_1,[~~]_1, \{~\}_1)$ and $(L_2,[~~]_2, \{~\}_2)$ respectively and $(V,m_l^1,m_r^1, m_l^2,m_r^2)$ and $(V,m_l^{1'},m_r^{1'}, m_l^{2'},m_r^{2'})$ be their corresponding representations given by Theorem \ref{Aex1}.\\
If $F$ is the isomorphism from $L_1$ to $L_2$ as given in the commutative diagram \ref{commdia} then 
define $s_1':L\to L_1$ by $s_1'=F^{-1} \circ s_2$. Then\\
$$\rho \circ s_1'=\rho \circ F^{-1} \circ s_2=\rho' \circ F \circ F^{-1} \circ s_2=\rho' \circ s_2=Id_L. $$
Thus $s_{1}'$ is a section of $L_1$ that gives the representation $(V,m_l^1,m_r^1, m_l^2,m_r^2)$ .\\
For all $x\in L, u\in V$ and $F|_V=Id_V$ we have
\begin{eqnarray}
   m_l^1(x,u)&=&[s_1'(x),u]_1\notag\\
   &=&[F^{-1} \circ s_2(x),u]_1\notag\\
   &=&F^{-1}[s_2(x),u]_2\notag\\
   &=&F^{-1}m_l^{1'}(x,u)=m_l^{1'}(x,u).\notag
\end{eqnarray}
    Similarly, we can show $m_r^1=m_r^{1'},~~m_l^2=m_l^{2'},~~m_r^2=m_r^{2'}$.


\thm With $\theta, \tilde{\theta},m_l^1,m_r^1, m_l^2,m_r^2$ defined above, $L\oplus V$ with the following bilinear operations is a compatible Leibniz algebra.
$$[x+u,y+v]_{\oplus}=[x,y]+\theta(x,y)+m_l^1(x,v)+m_r^1(u,y)$$
$$\{x+u,y+v\}_{\oplus}=\{x,y\}+\tilde \theta(x,y)+m_l^2(x,v)+m_r^2(u,y),$$
$~\forall x,y \in L,~u,v\in V.$
\proof The proof of the theorem is straightforward.

\thm $(\theta,\tilde \theta)$ is a $2$-cocycle of $(L,\pi_1,\pi_2)$ with coefficients in the representation $(V,m_l^1,m_r^1, m_l^2,m_r^2)$.
\proof For any $x,y,z \in L$
\begin{scriptsize}
\begin{eqnarray}
&&d^2_{\pi_1+m_l^1+m_r^1}\theta(x,y,z)\notag\\
&=&-[\hat \pi_1+\hat m_l^1+\hat m_r^1, \hat \theta]_B\notag \\
&=& m_r(\theta(x,y),z)-m_l(x,\theta(y,z))+m_l(y,\theta(x,z))+\theta ([x,y],z)-\theta(x,[y,z])+\theta(y,[x,z])\notag \\ 
&=&[\theta(x,y),s(z)]-[s(x),\theta(y,z)]+[s(y),\theta(x,z)]+[s[x,y],s(z)]-\notag\\
&&s[[x,y],z]-[s(x),s[y,z]]+s[x,[y,z]]+[s(y),s[x,z]]-s[y,[x,z]]\notag \\
&=&[[s(x),s(y)],s(z)]-\cancel{[s[x,y],s(z)]}-[s(x),[s(y),s(z)]]+\cancel{[s(x),s[y,z]]}+
[s(y),[s(x),s(z)]]-\notag \\
&&\cancel{[s(y),s[x,z]]}+\cancel{[s[x,y],s(z)]}-s[[x,y],z]-\cancel{[s(x),s[y,z]]}+s[x,[y,z]]+
\cancel{[s(y),s[x,z]]}-s[y,[x,z]]\notag \\
&=&[[s(x),s(y)],s(z)]-[s(x),[s(y),s(z)]]+[s(y),[s(x),s(z)]]-s[[x,y],z]+s[x,[y,z]]-s[y,[x,z]]\notag\\
&=&0.\notag
\end{eqnarray}
\end{scriptsize}
Similarly we can show $d^2_{\pi_2+m_l^2+m_r^2}\tilde \theta(x,y,z)=0.$
\begin{scriptsize}
\begin{eqnarray}
&&d^2_{\pi_2+m_l^2+m_r^2}\theta(x,y,z)+d^2_{\pi_1+m_l^1+m_r^1}\tilde \theta(x,y,z)\notag\\
&=&\{\theta(x,y),s(z)\}-\{s(x),\theta(y,z)\}+\{s(y),\theta(x,z)\}+\{s[x,y],s(z)\}-s[\{x,y\},z]-[s(x),s\{y,z\}]+\notag\\
&&s[x,\{y,z\}]+[s(y),s\{x,z\}]-s[y,\{x,z\}]+\notag\\
&&[\tilde \theta(x,y),s(z)]-[s(x),\tilde\theta(y,z)]+[s(y),\tilde\theta(x,z)]+\{s[x,y],s(z)\}-s\{[x,y],z\}-\{s(x),s[y,z]\}+\notag\\
&&s\{x,[y,z]\}+\{s(y),s[x,z]\}-s\{y,[x,z]\}\notag\\
&=&\{[s(x),s(y)],s(z)\}-\cancel{\{s[x,y],s(z)\}}-\{s(x),[s(y),s(z)]\}+\cancel{\{s(x),s[y,z]\}}+\{s(y),[s(x),s(z)]\}\notag\\
&&-\cancel{\{s(y),s[x,z]\}}+\cancel{\{s[x,y],s(z)\}}-s[\{x,y\},z]-\cancel{[s(x),s\{y,z\}]}+s[x,\{y,z\}]+\cancel{[s(y),s\{x,z\}]}\notag\\
&&-s[y,\{x,z\}]+[\{s(x),s(y)\},s(z)]-\cancel{[s\{x,y\},s(z)]}-[s(x),\{s(y),s(z)\}+\cancel{[s(x),s\{y,z\}]}\notag\\
&&+[s(y),\{s(x),s(z)\}]-\cancel{[s(y),s\{x,z\}]}+ \cancel{\{s[x,y],s(z)\}}-s\{[x,y],z\}-\cancel{\{s(x),s[y,z]\}}\notag\\
&&+s\{x,[y,z]\}+\cancel{\{s(y),s[x,z]\}}-s\{y,[x,z]\}\notag\\
&=&\{[s(x),s(y)],s(z)\}-\{s(x),[s(y),s(z)]\}+\{s(y),[s(x),s(z)]\}-s[\{x,y\},z]+s[x,\{y,z\}]\notag\\
&&-s[y,\{x,z\}]+[\{s(x),s(y)\},s(z)]-[s(x),\{s(y),s(z)\}+[s(y),\{s(x),s(z)\}]\notag\\
&&-s\{[x,y],z\}+s\{x,[y,z]\}-s\{y,[x,z]\}\notag\\
&=&0,~~~~~~~~~~~~~~~~\text{by the compatibility condition \eqref{CLA}}\notag
\end{eqnarray}
\end{scriptsize}
Thus $$\partial ^2(\theta,\tilde \theta)=(d^2_{\pi_1+m_l^1+m_r^1}\theta,d^2_{\pi_2+m_l^2+m_r^2}\theta+d^2_{\pi_1+m_l^1+m_r^1}\tilde \theta,d^2_{\pi_2+m_l^2+m_r^2}\tilde \theta )=0.$$

\lem \label{lem1}The cohomology class of the cocyle $(\theta,\tilde \theta)$ does not depend on the choice of sections.
\proof Let $(\theta,\tilde \theta)$ and $(\theta',\tilde \theta')$ be the cocycles corresponding to sections $s$ and $s'$.\\
Define 
$\psi: L\to V$ as $\psi(x)=s(x)-s'(x)$\\
Then 
\begin{eqnarray}
\theta(x,y)&=&[s(x),s(y)]_{\tilde L}-s[x,y]\notag\\
&=&[\psi(x)+s'(x), \psi(y)+s'(y)]_{\tilde L}-s'[x,y]-\psi[x,y]\notag\\
&=&[s'(x), \psi(y)]_{\tilde L}+[\psi (x), s'(y)]_{\tilde L}+[s'(x),s'(y)]_{\tilde L}-s'[x,y]-\psi[x,y]\notag\\
&=&m_l^1(x,\psi(y))+m_r^1(y,\psi (x))+\tilde \theta'(x,y)-\psi[x,y]\notag\\
&=&\partial_{\pi_1+m_l^1+m_r^1}\psi(x,y)+\tilde \theta'(x,y)\notag
\end{eqnarray}

Hence we have $\theta -\theta'=\partial_{\pi_1+m_l^1+m_r^2}\psi$.
Similarly it can be shown $\tilde \theta -\tilde \theta'=\partial_{\pi_2+m_l^2+m_r^2}\psi$.\\
Thus $(\theta -\theta',\tilde \theta -\tilde \theta')=\partial \psi$ i.e. $(\theta, \tilde \theta)$ and $(\theta ', \tilde \theta ')$ are in the same cohomology class.
\thm Equivalent abelian extensions give rise to the same cocycle $(\theta_1,\tilde \theta_1)$.
\proof Let $L_1$ and $L_2$ be equivalent abelian extensions of the compatible Leibniz algebra $(L,\pi_1, \pi_2)$ by $(V,0,0)$ as detailed by the commutative diagram
\begin{equation}\label{commdia}
\begin{tikzcd}
0\arrow{r} &  V \arrow{r}{i_1} \arrow{d}{id_v}  &  L_1 \arrow{r}{\rho_1} \arrow{d}{\phi} &  L \arrow{r}\arrow{d}{id_L} & 0\\
0\arrow{r} &  V \arrow{r}{i_2}                          &  L_2 \arrow{r}{\rho_2}                            &L\arrow{r}                         & 0.
\end{tikzcd}
\end{equation}
and sections $s_1:L\to L_1$ and $s_2:L\to L_2$ and cocycles $(\theta_1,\tilde \theta_1)$ and $(\theta_2,\tilde \theta_2)$ respectively.\\
We note that $$\rho_2(\phi \circ s_1)=(\rho_2 \circ \phi)s_1=\rho_1 \circ s_1=Id_L.$$
Hence $\phi\circ s_1$ is a section of $L_2$.\\
From Theorem \ref{lem1} we know that the cohomological class of the cocyle $(\theta,\tilde \theta)$ does not depend on the choice of sections. Hence we can take $s_2=\phi\circ s_1$.\\
Thus $\forall~x,y\in L$
\begin{eqnarray}
\theta_2(x,y)&=&[s_2(x),s_2(y)]_{2}-s_2[x,y]\notag\\
&=&[\phi\circ s_1(x),\phi \circ s_1(y)]_{2}-\phi \circ s_1[x,y]\notag\\
&=&\phi([ s_1(x), s_1(y)]_{1}-s_1[x,y])\notag\\
&=&[ s_1(x), s_1(y)]_{1}-s_1[x,y]\notag\\
&=&\theta_1(x,y).\notag
\end{eqnarray}
Likewise we can show $\tilde \theta_1=\tilde \theta_2.$
\thm Cohomologous cocycles $(\theta_1, \tilde \theta_1)$ and $(\theta_2, \tilde \theta_2)$ give rise to equivalent abelian extensions.
\proof Since $(\theta_1, \tilde \theta_1)$ and $(\theta_2, \tilde \theta_2)$ are cohomologous there exists $\phi: L \to V$ such that
$$\theta _1=\theta_2 +\partial_{\pi_1+m_l^1+m_r^1}\phi$$
$$\tilde \theta _1=\tilde \theta_2 +\partial_{\pi_2+m_l^2+m_r^2}\phi.$$
Suppose  $(L\oplus _1 V,[~~]_1,\{~\}_1)$ and $(L\oplus _2 V,[~~]_2,\{~\}_2)$ are abelian extensions of $(L,[~],\{~\})$ by $(V,0,0)$ with respect to $(\theta_1, \tilde \theta_1)$ and $(\theta_2, \tilde \theta_2)$  respectively.\\
 
Define $\psi: L\oplus_1 V \to L\oplus_2 V $ by $\psi(x+u)=x+u+\phi(x),~\forall x\in L, u\in V$.\\
For $x,y\in L,u,v \in V$,
\begin{eqnarray}
  &&\psi([x+u,y+v]_1)-[\psi(x+u),\psi (y+v)]_2\notag\\
  &=&\psi([x,y]+\theta_1(x,y)+m_l^1(x,v)+m_r^1(u,y))-[x+u+\phi(x), y+v+\phi(y)]_2\notag\\
   &=&[x,y]+\theta_1(x,y)+m_l^1(x,v)+m_r^1(u,y))+\phi[x,y]- \{[x,y]+\theta_2(x,y)+\notag\\
&&m_l^1(x,v)+m_l^1(x,\phi(y))+m_r^1(u,y)+m_r^1(\phi(x),y)\}\notag\\
   &=&\theta_1(x,y)+\phi[x,y]-\theta_2(x,y)-m_l^1(x,\phi(y))-m_r^1(\phi(x),y)\notag\\
   &=& \theta_1(x,y)-\theta_2(x,y)-\partial_{\pi_1+m_l^1+m_r^1}\phi(x,u)\notag\\
   &=&0.\notag
   \end{eqnarray}
   Similarly, we can show $\psi\{x+u,y+v\}_1=\{\psi(x+u),\psi (y+v)\}_2$.\\
   It is routine to verify that the following diagram commutes
   \begin{equation*}
   \begin{tikzcd}
0\arrow{r} &  V \arrow{r}{i_1} \arrow{d}{Id_v}  &  L\oplus_1 V \arrow{r}{\rho_1} \arrow{d}{\phi} &  L \arrow{r}\arrow{d}{Id_L} & 0\\
0\arrow{r} &  V \arrow{r}{i_2}                          & L\oplus_2 V  \arrow{r}{\rho_2}                            &L\arrow{r}                         & 0.
\end{tikzcd}
\end{equation*}


\begin{rem}
From the last two theorems we conclude that the 
 abelian extensions of compatible Leibniz algebra $(L,\pi_1,\pi_2)$ by $(V,0,0)$ are characterised by the $2^{nd}$ cohomology group $H^2(L,V)$.
\end{rem}



\end{document}